\documentclass[11pt]{article} 
\usepackage{amsfonts,amsmath,amsxtra}
\usepackage{latexsym}
\usepackage{amssymb}

\def\hybrid{\topmargin 0pt      \oddsidemargin 0pt
        \headheight 0pt \headsep 0pt
        \textwidth 16.5cm
        \textheight 23cm
        \marginparwidth 0.0in
        \parskip 5pt plus 1pt   \jot = 1.5ex}
\catcode`\@=11
\def\marginnote#1{}
\newcount\hour
\newcount\minute
\newtoks\amorpm
\hour=\time\divide\hour by60 \minute=\time{\multiply\hour by60
\global\advance\minute by-\hour}
\edef\standardtime{{\ifnum\hour<12 \global\amorpm={am}%
        \else\global\amorpm={pm}\advance\hour by-12 \fi
        \ifnum\hour=0 \hour=12 \fi
      \number\hour:\ifnum\minute<10 0\fi\number\minute\the\amorpm}}
\edef\militarytime{\number\hour:\ifnum\minute<10 0\fi\number\minute}

\def\draftlabel#1{{\@bsphack\if@filesw {\let\thepage\relax
   \xdef\@gtempa{\write\@auxout{\string
      \newlabel{#1}{{\@currentlabel}{\thepage}}}}}\@gtempa
   \if@nobreak \ifvmode\nobreak\fi\fi\fi\@esphack}
        \gdef\@eqnlabel{#1}}
\def\@eqnlabel{}
\def\@vacuum{}
\def\draftmarginnote#1{\marginpar{\raggedright\scriptsize\tt#1}}

\def\draft{\oddsidemargin -0.1truein
        \def\@oddfoot{\sl preliminary draft \hfil
        \rm\thepage\hfil\sl\today\quad\militarytime}
        \let\@evenfoot\@oddfoot \overfullrule 3pt
        \let\label=\draftlabel
        \let\marginnote=\draftmarginnote
\def\@eqnnum{{\rm (\theequation)}
\rlap{\kern\marginparsep\tt\@eqnlabel}%
\global\let\@eqnlabel\@vacuum}  }

\newfont{\Bbbb}{msbm7 scaled 1\@ptsize00}
\newcommand{\zs}{\raise-1pt\hbox{$\mbox{\Bbbb Z}$}}

scaled 1\@ptsize00

\font\sevenmsa=msam6 
\newfam\msafam
\textfont\msafam=\sevenmsa
\def\hexnumber@#1{\ifnum#1<10 \number#1\else
\ifnum#1=10 A\else\ifnum#1=11 B\else\ifnum#1=12 C\else \ifnum#1=13
D\else\ifnum#1=14 E\else\ifnum#1=15 F\fi\fi\fi\fi\fi\fi\fi}
\def\msa@{\hexnumber@\msafam}
\def\llcorner{\delimiter"4\msa@78\msa@78 }
\def\lrcorner{\delimiter"5\msa@79\msa@79 }
\mathchardef\blacktriangleright="3\msa@49
\mathchardef\blacktriangleleft="3\msa@4A \font\tenmsb=msbm10 scaled
1\@ptsize00
\newfam\msbfam
\textfont\msbfam=\tenmsb \scriptfont\msbfam=\tenmsb

\newdimen\Squaresize \Squaresize=14pt
\newdimen\Thickness \Thickness=0.5pt

\def\Square#1{\hbox{\vrule width \Thickness
   \vbox to \Squaresize{\hrule height \Thickness\vss
      \hbox to \Squaresize{\hss#1\hss}
   \vss\hrule height\Thickness}
\unskip\vrule width \Thickness} \kern-\Thickness}

\def\Vsquare#1{\vbox{\Square{$#1$}}\kern-\Thickness}

\def\numberbysection{\@addtoreset{equation}{section}
        \def\theequation{\thesection.\arabic{equation}}}
\numberbysection

\renewcommand{\theequation}{\thesection.\arabic{equation}}
\def\titlepage{\@restonecolfalse\if@twocolumn\@restonecoltrue\onecolumn
     \else \newpage \fi \thispagestyle{empty}\c@page\z@
        \def\thefootnote{\fnsymbol{footnote}} }

\def\endtitlepage{\if@restonecol\twocolumn \else  \fi
        \def\thefootnote{\arabic{footnote}}
        \setcounter{footnote}{0}}  
\relax

\hybrid
\makeatletter
\newdimen\normalarrayskip            
\newdimen\minarrayskip               
\normalarrayskip\baselineskip \minarrayskip\jot
\newif\ifold             \oldtrue            \def\new{\oldfalse}
\def\arraymode{\ifold\relax\else\displaystyle\fi}
\def\eqnumphantom{\phantom{(\theequation)}} 
\def\@arrayskip{\ifold\baselineskip\z@\lineskip\z@
     \else
     \baselineskip\minarrayskip\lineskip1\baselineskip\fi}

\def\@arrayclassz{\ifcase \@lastchclass \@acolampacol \or
\@ampacol \or \or \or \@addamp \or
   \@acolampacol \or \@firstampfalse \@acol \fi
\edef\@preamble{\@preamble
  \ifcase \@chnum
     \hfil$\relax\arraymode\@sharp$\hfil
     \or $\relax\arraymode\@sharp$\hfil
     \or \hfil$\relax\arraymode\@sharp$\fi}}

\def\@array[#1]#2{\setbox\@arstrutbox=\hbox{\vrule
     height\arraystretch \ht\strutbox
     depth\arraystretch \dp\strutbox
width\z@}\@mkpream{#2}\edef\@preamble{\halign \noexpand\@halignto
\bgroup \tabskip\z@ \@arstrut \@preamble \tabskip\z@ \cr}%
\let\@startpbox\@@startpbox \let\@endpbox\@@endpbox
  \if #1t\vtop \else \if#1b\vbox \else \vcenter \fi\fi
  \bgroup \let\par\relax
  \let\@sharp##\let\protect\relax
  \@arrayskip\@preamble}
\def\eqnarray{\stepcounter{equation}%
              \let\@currentlabel=\theequation
              \global\@eqnswtrue
              \global\@eqcnt\z@
              \tabskip\@centering              
              \let\\=\@eqncr
              $$%
            \halign to \displaywidth  \bgroup
             \eqnumphantom \@eqnsel
      \hskip\@centering                               
    $\displaystyle  \tabskip\z@ {##}$%
    &\global\@eqcnt\@ne \hskip 2\arraycolsep
         $ \displaystyle  \arraymode{##}$\hfil
    &\global\@eqcnt\tw@ \hskip 2\arraycolsep
         $\displaystyle\tabskip\z@{##}$\hfil
         \tabskip\@centering
    &{##}\tabskip\z@\cr}
\makeatother
\newcommand{\RR}{{\mathbb{R}}}

\def\IC{\mathbb{C}}

\def\IR{\mathbb{R}}

\def\CH {\mathcal{H}}

\def\CQ {\mathcal{Q}}

\def\CU {\mathcal{U}}
\def\CV {\mathcal{V}}

\def\CZ {\mathcal{Z}}

\def\la{\lambda}

\def\pr {\partial}

\def\Tr{{\rm Tr}}

\def\frak{\mathfrak}
\def\Fg{{\frak g}}

\newtheorem{te}{Theorem}[section]

\newtheorem{prop}{Proposition}[section]           

\newcommand\bqa{\begin{eqnarray}}
\newcommand\eqa{\end{eqnarray}}
\def\be{\begin{eqnarray}\new\begin{array}{cc}}
\def\ee{\end{array}\end{eqnarray}}

\def\beq{\begin{equation}}
\def\eeq{\end{equation}}
\def\bse{\begin{subequations}}                
\def\ese{\end{subequations}}
\def\bp{\begin{pmatrix}}
\def\ep{\end{pmatrix}}

\def\i{\imath}

\newcommand{\proof}{\noindent {\it Proof}. }
\def\stack#1#2{\raise0.7pt\hbox{$\mathrel{\mathop{#2}\limits^{#1}}$}}
\def\tr{\triangleright}
\def\tl{\triangleleft}
\def\sem{\mathsurround=0pt \raise1pt
\hbox{$\scriptscriptstyle>\!\!$}\:\!\!\tl}
\def\mes{\mathsurround=0pt \tr\!\:\!\raise0.8pt
\hbox{$\scriptscriptstyle\!\!<$}\,}
\def\]{\mathsurround=0pt ]\raise-2pt\hbox{$_\ast$}}

\def\<{\langle}
\def\>{\rangle}

\def\CQ{{\cal Q}}
\def\frak{\mathfrak}

\def\CU{{\cal U}}
\def\CZ{{\cal Z}}

\def\CH{\mathcal{H}}

\def\we{\raise-1pt\hbox{$\,\stackrel{\wedge}{,}\,$}}
\def\tr{{\rm tr}\,}
\def\Tr{{\rm Tr}\,}
\def\pr {\partial}

\newcounter{pac}[section]

\newcounter{pacc}[subsection]

0

0

\title{\bf On universal  Baxter operator for classical groups}
\begin{document}
\author{Anton A. Gerasimov and  Dimitri R. Lebedev}
\date{}

\maketitle

\renewcommand{\abstractname}{}

\begin{abstract}

\noindent {\bf Abstract}.
The universal Baxter operator  is an element
 of the Archimedean spherical Hecke algebra $\CH(G,K)$, $K\subset
 G$ be a  maximal compact subgroup of a Lie group $G$. It
 has a defining property to act in
spherical principle series representations of $G$ via multiplication
on the corresponding  local Archimedean $L$-factors.
Recently such operators were introduced for $G=GL_{\ell+1}(\IR)$
as generalizations of the Baxter operators arising in the theory of quantum
Toda chains. In this note
we provide universal Baxter operators for
classical groups $SO_{2\ell}$, $Sp_{2\ell}$ using the results of
Piatetski-Shapiro and Rallis on integral representations of
local Archimedean $L$-factors.

\end{abstract}
\vspace{5 mm}

\maketitle

\renewcommand{\abstractname}{}

\section{Introduction}

Interactions between theory of quantum integrable systems
and representation theory were intensive and mutually fruitful
starting from the earlier studies of quantum integrable systems.
Interpretation of the main technical tools used to solve quantum
integrable systems in terms of representation theory
was always a challenge promising  not only  a better understanding of
the phenomena of quantum integrability  but hopefully providing new methods
in representation theory as well. 
One of the important results in the theory of quantum integrable systems
is a construction of a class of operators introduced by Baxter \cite{B}.
The fundamental  role playing by the Baxter operators in
finding explicit solutions of quantum integrable systems
becomes more and more obvious.
However  a representation theory interpretation of
these operators is  not quite satisfactory. Recently a progress in this
direction was achieved for the Baxter operators associated
with the  quantum $\mathfrak{gl}_{\ell+1}$-Toda chains \cite{GLO2}.
Recall that the ring of    $\mathfrak{g}$-Toda chain quantum  Hamiltonians
can be identified with the center $\CZ$ of the universal enveloping
algebra $\CU \Fg$. Common eigenfunctions of the quantum Hamiltonians
are given by the $\mathfrak{g}$-Whittaker
functions identified with particular matrix elements of
principle series  representations  of $\CU\Fg$.
In \cite{GLO2} the integral Baxter  operator  for $\mathfrak{gl}_{\ell+1}$-Toda chain
was introduced (see \cite{PG} for analogous  Baxter
operators for affine $\widehat{\mathfrak{gl}}_{\ell+1}$-Toda
chains). The integral Baxter operator
is a one-parameter family of the integral operators
such that the $\mathfrak{gl}_{\ell+1}$-Whittaker
functions are their common eigenfunctions. Taking into
account the representation theory interpretation of the quantum Hamiltonians
as elements of the center $\CZ\subset \CU \Fg$ it is
natural  to look  for a similar  representation theory
interpretation of the Baxter operators for quantum
$\mathfrak{gl}_{\ell+1}$-Toda chains.
It was argued in \cite{GLO2} that the natural representation theory framework for
the Baxter operators is a
theory of the spherical Hecke algebras
$\CH=\CH(GL_{\ell+1}(\IR),K)$ where $K$ is a maximal compact subgroup $K$ of $G(\IR)$.
Note that the Whittaker functions
can be understood as functions on $G(\IR)$ via matrix element
interpretation mentioned above.
Then  the integral Baxter operator of
$\mathfrak{gl}_{\ell+1}$-Toda chain arises via a convolution of
an explicitly defined one-parameter family of elements of $\CH$
with the $\mathfrak{gl}_{\ell+1}$-Whittaker functions. The kernels of the
 corresponding integral operators are given basically by
the Gaussian $K$-biinvariant measures on $GL_{\ell+1}(\IR)$.
Such family of elements  of the spherical Hecke algebra $\CH$ was
called the universal
Baxter operator. A reason for this term is that
such operator provides simultaneously a construction of
the Baxter operators for a whole class of integrable systems
associated with a pair $K\subset GL_{\ell+1}(\IR)$
(for example including along with $\mathfrak{gl}_{\ell+1}$-Toda chains also
$\mathfrak{gl}_{\ell+1}$-Calogero-Sutherland systems).
The surprising result of \cite{GLO2}
is that the eigenvalues of the universal Baxter operators acting
on the $\mathfrak{gl}_{\ell+1}$-Whittaker
functions corresponding to principle series representations $\CV$ of
$GL_{\ell+1}(\IR)$ are equal to local Archimedean $L$-factor attached
to $\CV$ via the local Archimedean Langlands correspondence.
Such characterization of the local Archimedean $L$-factor is
new and is in perfect correspondence with  the analogous constructions over
non-Archimedean fields.

In this short note we propose a  construction of
universal Baxter operator for classical series $Sp_{2\ell}$ and
$SO_{2\ell}$ generalizing results of \cite{GLO2}.
The construction of the universal Baxter operators (Theorem \ref{Main})
essentially  relies on the results of \cite{PSR}  on
integral representations for local Archimedean $L$-factors
associated with principle series representations of  $SO_{2\ell}$ and
$Sp_{2\ell}$ (see \cite{LR} for an extension to all classical series). 
We also show that some of the results of \cite{PSR} can be easily
rederived using the approach put forward in \cite{GLO2}.
Let us stress that
 the case of the classical groups other then $GL_{\ell+1}$
reveals a new phenomena.
In \cite{GLO2} we argue that the Baxter integral operators  are a close cousins of
the recursive operators for the Whittaker functions. Recursive operators
for all classical series of finite Lie groups were constructed in
\cite{GLO1} and for classical groups other then $GL_{\ell+1}$
their integral  kernels are given by non-trivial integral expressions.
One should expect that a similar phenomena takes place  for the
universal Baxter  elements of the Hecke algebra $\CH(G,K)$ for $G$ being
classical groups other then $GL_{\ell+1}$.
An explicit form of the integral kernel given in Theorem \ref{Main} confirm these
expectations. Detailed discussion  of the proposed construction
of the universal Baxter operators for all classical Lie groups
including a derivation of explicit
integral expressions for the Baxter operators  acting on
Whittaker functions  will appear elsewhere \cite{GLO3}.

{\em Acknowledgments}: The research was supported by  grant
RFBR-09-01-93108-NCNIL-a. The research of AG was  also partly
supported by Science Foundation Ireland grant.

\section{ Baxter operator and spherical Hecke algebras}

In this Section we review two particular classes of integral
representations of $\mathfrak{gl}_{\ell+1}$-Whittaker functions.
Let $E_{ij}$, $i,j=1,\ldots \ell+1$ be the standard basis of the Lie
algebra $\mathfrak{gl}_{\ell+1}$.
Let $\CZ(\CU\mathfrak{gl}_{\ell+1})\subset
\CU\mathfrak{gl}_{\ell+1}$ be a center of the  universal
enveloping algebra $\CU\mathfrak{gl}_{\ell+1}$.
Let  $B_{\pm}\subset GL_{\ell+1}(\IC)$
be  upper-triangular and lower-triangular
Borel subgroups and  $N_{\pm}\subset B_{\pm}$
be  upper-triangular and lower-triangular
unipotent subgroups. Denote by $\mathfrak{b}_{\pm}={\rm Lie}(B_{\pm})$ and
$\mathfrak{n}_{\pm}={\rm Lie}(N_{\pm})$ their Lie algebras.
 Let $\mathfrak{h}\subset \mathfrak{gl}_{\ell+1}$
be a diagonal Cartan subalgebra and $W=\mathfrak{S}_{\ell+1}$ be the Weyl
group of $GL_{\ell+1}$.   Using the Harish-Chandra
isomorphism of $\CZ(\CU\mathfrak{gl}_{\ell+1})$ with
$\mathfrak{S}_{\ell+1}$-invariant subalgebra
of the symmetric algebra $S^*\mathfrak{h}$  we  identify central characters with
homomorphisms $c:\IC[h_1,\cdots, h_{\ell+1}]^{\mathfrak{S}_{\ell+1}}\to \IC$.
Let $\pi_{\underline{\lambda}}:\CU\mathfrak{gl}_{\ell+1}\to
{\rm End}(\CV_{\underline{\la}})$,
$\CV_{\underline{\la}}={\rm   Ind}_{\CU\frak{b}_-}^{\CU
  \frak{gl}_{\ell+1}}
\chi_{\underline{\lambda}}$
be a family of principal series representations   of
$\CU\frak{gl}_{\ell+1}$  induced from  one-dimensional
representations
$\chi_{\underline{\lambda}}(b)=\prod_{j=1}^{\ell+1}\,|b_{jj}|^{
\imath \lambda_k-\rho_k}$
of $\CU\frak{b}_-$. Here
$\underline{\lambda}=(\lambda_1,\ldots, \lambda_{\ell+1})\in
\IR^{\ell+1}$  and $\rho_k=(\ell-2k+2)/2,\,\,\,\,
k=1,\ldots,\ell+1$.  Denote $\<\,,\,\>$ a pairing
 on $\CV_{\underline{\lambda}}$.
We suppose that the action of the Cartan subalgebra $\mathfrak{h}$
 in representation  $\CV_{\underline{\lambda}}$ can be integrated
to an  action of the corresponding Cartan subgroup $H\subset GL_{\ell+1}(\IC)$.

The $\mathfrak{gl}_{\ell+1}$-Whittaker function can be defined
as  a matrix element of a principle series representation
$\CV_{\underline{\lambda}}$ of $G=GL_{\ell+1}(\IR)$ (see e.g. \cite{J}, \cite{Ha}).
Let us fix the Iwasawa decomposition $G=N_-AK-$ where $K$ is a maximal
compact subgroup of  $G$, $A$ be a group of diagonal matrices with the
positive elements and $N_-$ is a maximal unipotent subgroup of
lower-triangular matrices with unites on the diagonal.
Let $\rho(g)$ be   given by $\rho(nak)=\<\rho,\log a\>$
 where $k\in K$, $a\in A$ and $n\in N_-$. Let  $\phi_K$
be  a spherical vector in $\mathcal{V}_{\underline{\la}}$ and $\psi_{N_-}$
be the  Whittaker vector defined by the condition
\be
\psi_{N_-}(bgn)=\chi_{\underline{\la}}(b)\chi_{N_-}(n)\,\psi_{N_-}(g),\qquad
 n\in N_-,\quad b\in B_-,
\ee
\be
\phi_K(bgk)=\chi_{\underline{\la}}(b)\phi_K(g),\qquad \qquad
k\in K,\,b\in B_-,
\ee
where $\chi_{N_-}(n)=\exp(2\pi\imath\sum_{j=1}^{\ell}n_{j+1,j})$.
Then the Whittaker function on $G$ is defined as a matrix
element
\be\label{Wf}
\Psi^{\mathfrak{gl}_{\ell+1}}_{\underline{\la}}(g)
= e^{\rho(g)}\<\psi_{N_-},\pi_{\underline{\la}}(g) \,\phi_K\>.
\ee
The function \eqref{Wf} satisfies the obvious functional equation
\bqa \Psi^{\mathfrak{gl}_{\ell+1}}_{\underline{\la}}(ngk)
=\chi_{N_-}(n)\,
\Psi^{\mathfrak{gl}_{\ell+1}}_{\underline{\la}}(g),\qquad g\in G,\,
k\in K,\,n\in N_-,
\eqa
and  descends to a function  on the space $A$ of the
diagonal matrices $a={\rm diag}(e^{x_1},\ldots,
e^{x_{\ell+1}})$.

Standard considerations (see e.g. \cite{STS})
show that  the matrix element \eqref{Wf}  is a common
eigenfunction  of a family of commuting differential operators
descending  from   generators of
$\CZ(\CU\mathfrak{gl}_{\ell+1})$ acting in $\CV_{\underline{\lambda}}$.
These differential operators  can be identified with quantum Hamiltonians of
$\mathfrak{gl}_{\ell+1}$-Toda chain. For example the simplest non-trivial quantum
Hamiltonian acts on the Whittaker function via the differential operator
$$
H_2^{\mathfrak{gl}_{\ell+1}}=-\frac{1}{2}\sum_{i=1}^{\ell+1}\frac{\pr^2}{\pr
  x_i^2}+4\pi^2\sum_{i=1}^{\ell}e^{2(x_{i+1}-x_{i})}.
$$
Using explicit realizations of the
universal enveloping algebra representation
$\pi_{\underline{\lambda}}$ via difference/
differential operators acting in an appropriate
space of functions the matrix element representation  \eqref{Wf}
leads to various  integral representations of
the Whittaker function.

As it was demonstrated in \cite{GLO2}
the $\mathfrak{gl}_{\ell+1}$-Whittaker function being a  common
eigenfunction of a family of mutually commuting differential operators
is also a common eigenfunction of a one-parameter family of
integral operators. These integral
operators were  called the Baxter operators due to their relation with the
Baxter operators in the theory of quantum integrable systems.

Define the Baxter $Q$-operator
as an integral operator
acting in an appropriate space of functions of $\ell+1$ variables
with the integral  kernel
\be\label{Baxter} Q^{\mathfrak{gl}_{\ell+1}}
(\underline{x},\underline{y}|s)=2^{\ell+1}\,\exp\Big\{\sum_{j=1}^{\ell+1}(\i
s -\rho_j) (x_j-y_j)-\\  -\pi
\sum_{k=1}^{\ell}\Big(e^{2(x_k-y_k)}+e^{2(y_{k+1}-x_{k})}\Big)- \pi
e^{2(x_{\ell+1}-y_{\ell+1})}\Big\},
\ee
 where $\rho=(\rho_1,\ldots, \rho_{\ell+1})\in \RR^{\ell+1}$ , with
 $\rho_j=\frac{\ell}{2}+1-j,\,\,\,j=1,\ldots,\ell+1$.

The operator  $Q^{\mathfrak{gl}_{\ell+1}}(s)$
satisfies the following commutativity relations:
\bqa\label{firstpr}
Q^{\mathfrak{gl}_{\ell+1}}(s)\cdot
Q^{\mathfrak{gl}_{\ell+1}}(s')=
Q^{\mathfrak{gl}_{\ell+1}}(s')\cdot
Q^{\mathfrak{gl}_{\ell+1}} (s), \eqa
\bqa\label{secondpr}
Q^{\mathfrak{gl}_{\ell+1}}(s)\cdot
H_r^{\mathfrak{gl}_{\ell+1}}=
H_r^{\mathfrak{gl}_{\ell+1}}\cdot
Q^{\mathfrak{gl}_{\ell+1}}(s),\qquad r=1,\ldots \ell+1,
\eqa
where $H_r^{\mathfrak{gl}_{\ell+1}}$ are quantum Hamiltonians of
$\mathfrak{gl}_{\ell+1}$-Toda chain.

With respect to $Q^{\mathfrak{gl}_{\ell+1}}(s)$
the $\mathfrak{gl}_{\ell+1}$-Whittaker function \eqref{Wf}
has  the following eigenfunction identity:
\bqa\label{eigenprop}
Q^{\mathfrak{gl}_{\ell+1}}(s)\,
\cdot \Phi^{\mathfrak{gl}_{\ell+1}}_{\underline{\lambda}}
(\underline{x}):=
\int_{\RR^{\ell+1}}\,\prod_{i=1}^{\ell+1}\,dy_{i}\,\,
Q^{\mathfrak{gl}_{\ell+1}}(\underline{x},
\,\underline{y}|\,s)\,
\Phi^{\mathfrak{gl}_{\ell+1}}_{\underline{\lambda}}(\underline{y})\,=\,
L(s,\underline{\lambda})\,\,
\Phi^{\mathfrak{gl}_{\ell+1}}_{\underline{\lambda}}(\underline{x}),
\eqa
where
 \be\label{newW}
\Phi^{\mathfrak{gl}_{\ell+1}}_{\underline{\lambda}}(\underline{x})=
e^{-\<\rho,x\>}\,
\Psi^{\mathfrak{gl}_{\ell+1}}_{\underline{\lambda}}(\underline{x}).
\ee
We use the following notations $\underline{x}=(x_1,\ldots, x_{\ell+1})$,
$\underline{y}=(y_1,\ldots, y_{\ell+1})$,
$\underline{\la}=(\la_1,\ldots, \la_{\ell+1})$
and the eigenvalue in \eqref{eigenprop} is  given by
\be\label{lAL}
L(s,\underline{\lambda})=\prod_{j=1}^{\ell+1}\,
\pi^{-\frac{\imath  s -\imath  \lambda_j}{2}}
\Gamma\Big(\frac{\imath  s- \imath \lambda_{j}}{2}\Big).
\ee
The eigenvalue \eqref{lAL} is
a  local Archimedean $L$-factor associated with principle series
representation $\CV_{\underline{\lambda}}$.

The appearance of the local Archimedean  $L$-factors in
\eqref{eigenprop} is not accidental and is
related with the fact that  the integral operators \eqref{Baxter}
are realizations of particular  elements of
 the local spherical Archimedean Hecke algebra (see \cite{GLO2}
for detailed discussion). Recall that the  Archimedean Hecke algebra
$\CH=\CH(G(\mathbb{R}),K)$,  $K$ being  a maximal compact
subgroup of  $G(\RR)$   is defined
as an algebra of $K$-biinvariant functions
 on $G=G(\IR)$,  $\phi(g)=\phi(k_1gk_2)$, $k_1,k_2\in K$  acting by a convolution
\be\label{con5}
\phi*f(g)=\int_G \phi(g\tilde{g}^{-1})\,f(\tilde{g})
d\tilde{g}=\int_G\,\phi(\tilde{g})\,f(\tilde{g}^{-1}g)\,d\tilde{g},
\ee
where the last equality holds for unimodular groups. 
To ensure the convergence of
the integrals \eqref{con5} we imply that the elements of spherical
Hecke algebra belong to the Schwarz functional space
i.e. the functions such all their derivatives decrease at infinity
faster then inverse power of any polynomial. Spherical Hecke algebra
is isomorphic to the algebra
of $Ad_{G^{\vee}}$-invariant functions on $\mathfrak{g}^{\vee}={\rm
  Lie}(G^{\vee})$ where $G^{\vee}$ is  a  complex Lie group   dual to $G$
(e.g.  $A_{\ell}$, $B_{\ell}$,  $C_{\ell}$,  $D_{\ell}$ are dual to
 $A_{\ell}$, $C_{\ell}$,  $B_{\ell}$,  $D_{\ell}$ respectively).
Given a finite-dimensional representation
$\rho_V: G^{\vee}\to GL(V,\IC)$ one attaches
a local Archimedean  $L$-function  corresponding to a
 spherical irreducible representation of $G$ as follows (see e.g. \cite{Bu}, \cite{L}):
\bqa\label{archL}
L(s',\phi,\rho_V)=\prod_{j=1}^{\ell+1}
\pi^{-\frac{s'-\lambda'_j}{2}}\,
\Gamma\Big(\frac{s'-\lambda'_j}{2}\Big)
=\det_V\Big(\pi^{-\frac{s'-\rho_V(t_{\infty})}{2}}\,
\Gamma\Big(\frac{s'-\rho_V(t_{\infty})}{2}\Big)\Big),
\eqa
where $\rho_V(t_{\infty})={\rm diag}(\lambda_1',\ldots \lambda_{\ell+1}')$
and  $t_{\infty}$ is a conjugacy class in the Lie
 algebra $\mathfrak{g}^{\vee}={\rm Lie}(G^{\vee})$. In the case of
 $G=GL_{\ell+1}$, $V=\IC^{\ell+1}$, $s'=\imath s$, $\lambda_j'=\imath
 \lambda_j$   we recover \eqref{lAL}.

By the multiplicity one theorem \cite{Sha}, there is a unique
smooth  $K$-spherical vector $\phi_K$ in a principal series irreducible
representation $\mathcal{V}_{\underline{\lambda}}={\rm
Ind}_{B_-}^G\,\chi_{\underline{\lambda}}$. The action of a
$K$-biinvariant function $\phi$ on the spherical vector $\phi_K$ in
$\mathcal{V}_{\underline{\lambda}}$ is given by  the multiplication
by  a character $\Lambda_{\phi}$ of the Hecke algebra $\CH$:
\be\label{conS}
 \phi*\phi_K(g)=\int_G\,\, dg_1\,\,
\phi(g g_1^{-1})\,\phi_K(g_1)=\Lambda_{\phi}(\underline{\lambda}) \phi_K(g).
\ee
In particular, the  elements $\phi$ of the Hecke algebra act
via convolution on the Whittaker function \eqref{newW}  as follows:
\bqa\label{phiEigen}
 \Phi^{\mathfrak{gl}_{\ell+1}}_{\underline{\lambda}}
*\phi(g)=\Lambda_{\phi}(\underline{\lambda})\,\,
 \Phi^{\mathfrak{gl}_{\ell+1}}_{\underline{\lambda}}
(g),\qquad \phi\in \CH.
\eqa
Here the Whittaker function
$\Phi_{\underline{\lambda}}^{\mathfrak{gl}_{\ell+1}}$ is considered
as a function on $G$ such that \be\label{equivar}
\Phi^{\mathfrak{gl}_{\ell+1}}_{\underline{\lambda}}(nak)=\chi_{N_-}(n)\,
 \Phi^{\mathfrak{gl}_{\ell+1}}_{\underline{\lambda}}(a),
\ee where $\chi_{N_-}(n)$ is a character of $N_-$ trivial on $[N_-,N_-]$ and
$nak\in N_-AK$ is the Iwasawa decomposition of $G$.

One can express the eigenvalue $\Lambda_{\phi}(\underline{\lambda})$
corresponding to the action \eqref{phiEigen}  of an arbitrary element $\phi\in \CH$
as follows.   Consider an action of $\phi$
on a normalized spherical function $\varphi_{\underline{\lambda}}(g)$
 in a principle series representation
$\CV_{\underline{\lambda}}$ i.e.  $\varphi_{\underline{\lambda}}$ is a matrix element in
$\CV_{\underline{\lambda}}$   such that
$$
\varphi_{\underline{\lambda}}(k_1\,g \,
k_2)=\varphi_{\underline{\lambda}}(g),\qquad
\varphi_{\underline{\lambda}}(e)=1,
 \qquad g\in G, \quad k_1,k_2\in K.
$$
An explicit integral
representation for the normalized spherical function $\varphi_{\underline{\lambda}}(g)$ is
given by
\be
\varphi_{\underline{\lambda}}(g)=\int_K \,dk\,\,
e^{\<h(kg),\i\underline{\lambda}-\rho\>},
\ee
where $h(g)=\log a$,  $g=nak\in N_-AK$ is
the Iwasawa decomposition  of $G$ and we normalize the volume of the compact
subgroup as $\int_K dk=1$.
For an eigenvalue $\Lambda_{\phi}$ of $\phi$
$$
 \,\varphi_{\underline{\lambda}}*\phi(g)=\,
\Lambda_{\phi}(\underline{\lambda})\,\,\varphi_{\underline{\lambda}}(g),
$$
we then obtain  the integral representation
\be\label{Lfactor}
\Lambda_{\phi}(\underline{\lambda})= \,\varphi_{\underline{\lambda}}*\phi(e)=
\int_G dg
_{}\phi(g^{-1})\varphi_{\underline{\lambda}}(g).
\ee

In \cite{GLO2} it was shown
that  the local Archimedean $L$-factor can be
understand as an eigenvalue of a particular one-parameter
family of elements $\CQ^{\mathfrak{gl}_{\ell+1}}(s)$
of $\CH(GL_{\ell+1}(\IR),O_{\ell+1})$ acting on
$\mathfrak{gl}_{\ell+1}$-Whittaker functions
\be
 \Phi^{\mathfrak{gl}_{\ell+1}}_{\underline{\lambda}}*\CQ^{\mathfrak{gl}_{\ell+1}}(s)
 (g)=
\prod_{j=1}^{\ell+1}\,
\pi^{-\frac{\imath s-\imath\lambda_j}{2}}\,
\Gamma\Big(\frac{\imath s-\imath \lambda_{j}}{2}\Big)
\,\,\Phi^{\mathfrak{gl}_{\ell+1}}_{\underline{\lambda}}(g), \ee and
the action $\phi$ to the subspace of
functions satisfying (\ref{equivar}) coincides with  the action of the
integral operator $\CQ^{\mathfrak{gl}_{\ell+1}}(s)$ given explicitly by  \eqref{Baxter}.
The corresponding
one-parameter   family was called universal Baxter operator.

\begin{te}[GLO]\label{MainTh1}
 Let $\CQ^{\mathfrak{gl}_{\ell+1}}(s)$ be a $K$-biinvariant
function on $G=GL_{\ell+1}$  given by
\be\label{UBO}
\CQ^{\mathfrak{gl}_{\ell+1}}(g,s)=2^{\ell+1}\,|\det g|^{\imath
s+\frac{\ell}{2}} e^{-\pi{\rm Tr} g^tg}.
\ee
Then the action of $\CQ^{\mathfrak{gl}_{\ell+1}}(s)$ on the
Whittaker function $\Phi^{\mathfrak{gl}_{\ell+1}}_{\underline{\lambda}}(g)$ by a
convolution  descends to the action of
$Q^{\mathfrak{gl}_{\ell+1}}(s)$ with the integral kernel
\eqref{Baxter} and satisfies the relation
\be\label{Eigenprop}
 \Phi^{\mathfrak{gl}_{\ell+1}}_{\underline{\lambda}}*\CQ^{\mathfrak{gl}_{\ell+1}}(s)(g)=
L(s,\underline{\lambda})\,\,\Phi^{\mathfrak{gl}_{\ell+1}}_{\underline{\lambda}}(g),
\ee where $L(s,\underline{\lambda})$ is the  local Archimedean $L$-factor
\be\label{Lfac} L(s,\underline{\lambda})=\prod_{j=1}^{\ell+1}\,
\pi^{-\frac{\imath s -\imath \lambda_j}{2}}
\Gamma\Big(\frac{\imath s-\imath \lambda_{j}}{2}\Big).
\ee
\end{te}
Below we prove a slightly modified version of this Theorem.

\begin{te}\label{MainTh2}
 Let $\tilde{\CQ}^{\mathfrak{gl}_{\ell+1}}(s)$ be a $K$-biinvariant
function on $G=GL_{\ell+1}$  given by
\be\label{UBO1}
\widetilde{\CQ}^{\mathfrak{gl}_{\ell+1}}(g,s):=
\CQ^{\mathfrak{gl}_{\ell+1}}(g^{-1},s)=2^{\ell+1}\,
|\det g^{-1}|^{\imath
 s+\frac{\ell}{2}} e^{-\pi{\rm Tr} (g^{-1})^tg^{-1}}.
\ee
Then the  action of $\tilde{\CQ}^{\mathfrak{gl}_{\ell+1}}(g,s)$ on the  modified
Whittaker function $\Phi^{\mathfrak{gl}_{\ell+1}}_{\underline{\lambda}}(g)$ by a
convolution  descends to the action of
$\tilde{Q}^{\mathfrak{gl}_{\ell+1}}(s)$ with the integral kernel
\be\label{twQ}
\tilde{Q}^{\mathfrak{gl}_{\ell+1}}
(\underline{x},\underline{y}|s)=2^{\ell+1}\,\exp\Big\{\sum_{j=1}^{\ell+1}(\i
s +\rho_j) (y_j-x_j)-\\ \nonumber -\pi
\sum_{k=1}^{\ell}\Big(e^{2(y_k-x_k)}+e^{2(x_{k+1}-y_{k})}\Big)- \pi
e^{2(y_{\ell+1}-x_{\ell+1})}\Big\}, \ee
  and satisfies the relation
\be\label{EigenpropT}
\Phi^{\mathfrak{gl}_{\ell+1}}_{\underline{\lambda}}*\widetilde{\CQ}^{\mathfrak{gl}
_{\ell+1}}(s)(g)=
L(s,-\underline{\lambda})\,\,
\Phi^{\mathfrak{gl}_{\ell+1}}_{\underline{\lambda}}(g),
\ee where $L(s,\underline{\lambda})$ is the  local Archimedean
$L$-factor given by \eqref{Lfac}.
\end{te}

\noindent {\it Proof}.  The convolution of a $K$-biinvariant
function with  $\mathfrak{gl}_{\ell+1}$-Whittaker functions is given by
\be

\Phi^{\mathfrak{gl}_{\ell+1}}_{\underline{\lambda}}*\widetilde{\CQ}^{\mathfrak{gl}
_{\ell+1}}(g)
=\int_G\, d\tilde{g}\, \tilde{\CQ}^{\mathfrak{gl}_{\ell+1}}(g\tilde{g}^{-1})\,
\Phi^{\mathfrak{gl}_{\ell+1}}_{\underline{\lambda}}(\tilde{g})=\int_G\,
 d\tilde{g}\,
\tilde{\CQ}^{\mathfrak{gl}_{\ell+1}}
(g\tilde{g}^{-1})\,\<\psi_{N_-},\pi_{\underline{\lambda}}(\tilde{g})\phi_K\>.
\ee
Fix the Iwasawa decomposition
$\tilde{g}=\tilde{n}\tilde{a}\tilde{k}$, $\tilde{k}\in K$,
$\tilde{a}\in A$, $\tilde{n}\in N_-$ of a generic element
$\tilde{g}\in G$  and let
$\delta_{B_-}(\tilde{a})=\det_{\mathfrak{n}_-} {\rm
Ad}_{\tilde{a}}$.
We normalize the volume of the compact subgroup as $\int_K dk=1$.
We shall use  the notation
$d^{\times}a=da\cdot\det(a)^{-1}$ for $a\in A$.  Then for $a\in A$ we
have
\be
\Phi^{\mathfrak{gl}_{\ell+1}}_{\underline{\lambda}}*
\tilde{\CQ}^{\mathfrak{gl}_{\ell+1}}(a)=
\int_{AN_-}\, d^{\times}\tilde{a}d\tilde{n}\,\delta_{B_-}(\tilde{a})
\,\tilde{\CQ}^{\mathfrak{gl}_{\ell+1}}(a\tilde{n}^{-1}\tilde{a}^{-1})\,
\chi_{N_-}(\tilde{n})\,\Phi^{\mathfrak{gl}_{\ell+1}}_{\underline{\lambda}}
(\tilde{a})\,=
\\
=\int_{A}d^{\times}\tilde{a}\,\, K_{\tilde{\phi}}(a,\tilde{a})\,
\Phi^{\mathfrak{gl}_{\ell+1}}_{\underline{\lambda}} (\tilde{a}),\ee
with \be
K_{\tilde{\phi}}(a,\tilde{a})=\int_{N_-}\!\!d\tilde{n}\,\,\delta_{B_-}(\tilde{a})\,
\,{\CQ}^{\mathfrak{gl}_{\ell+1}}
(\tilde{a}\tilde{n}{a}^{-1})\,\chi_{N_-}(\tilde{n}), \\ \nonumber
\chi_{N_-}(\tilde{n})=\exp\Big\{\,2\pi \i
\sum_{i=1}^{\ell}\tilde{n}_{i+1,i}\Big\}.\ee  Thus to prove the
expression \eqref{twQ} for the integral kernel we should prove the following
\be
\tilde{Q}^{\mathfrak{gl}_{\ell+1}}(\underline{x},\underline{y}|s)
=\int_{N_-}\!\!d\tilde{n}\,\,\delta_{B_-}(\tilde{a})\,
{\CQ}(\tilde{a}\tilde{n}{a}^{-1}|s)\,\chi_{N_-}(\tilde{n}), \ee
where
\be
a={\rm\,diag}\,(e^{x_1},\ldots, e^{x_{\ell+1}}), \qquad
\tilde{a}={\rm diag}(e^{y_1},\ldots ,e^{y_{\ell+1}}), \qquad\\
\delta_{B_-}(\tilde{a})=e^{2\<\rho,\log
\tilde{a}\>}=e^{\sum_{i<j}(y_i-y_j)}. \ee
 For
$g=\tilde{a}\tilde{n}{a}^{-1}$ we have \be
\det\,g=e^{\sum_{i=1}^{\ell+1}(y_i-x_i)}, \qquad
\Tr\,g^tg=\sum_{i=1}^{\ell+1}\,e^{2(y_i-x_i)}+ \sum_{i<j}
\tilde{n}_{ij}^2e^{2(y_i-x_j)}, \ee where $\tilde{n}\in N_-$.
Taking into account that $\chi_{N_-}(\tilde{n})=\exp(
2\pi \i\sum_{i=1}^{\ell}\tilde{n}_{i+1,i})$ we obtain
\bqa
\tilde{Q}^{\mathfrak{gl}_{\ell+1}}(\underline{x},\underline{y}|s)=2^{\ell+1}\,
\int_{N_-}\!\! d\tilde{n}\,\,
e^{\sum_{i<j}(y_i-y_j)} \,e^{2\pi\i\sum_{k=1}^{\ell}\tilde{n}_{i+1,i}}\times \\
\nonumber\exp\Big\{\sum_{i=1}^{\ell+1}(\imath s+\frac{\ell}{2})
(y_i-x_i)- \pi\sum_{i=1}^{\ell+1}e^{2(y_i-x_i)}-
\pi\sum_{i<j}\tilde{n}_{ij}^2e^{2(y_i-x_j)}\Big\}=\\
2^{\ell+1}\,\exp\Big\{(\imath s+\frac{\ell}{2})\sum_{i=1}^{\ell+1}(y_i-x_i)-
\pi
\sum_{i=1}^{\ell+1}e^{2(y_i-x_i)}\Big\}\,e^{\sum_{i<j}(y_i-y_j)}\times \\
\nonumber \int_{\mathbb{R}^\ell}\prod_{i=1}^\ell\!d\tilde{n}_{i+1,i}\,\,
\exp\Big\{2\pi \i\sum_{k=1}^{\ell}\tilde{n}_{i+1,i}-\pi\sum_{i=1}^\ell
\tilde{n}_{i+1,i}^2e^{2(y_{i}-x_{i+1})}\Big\}\times \\ \nonumber
\prod_{i>j+1}\int\!\!d\tilde{n}_{ij}\,\, \exp\Big\{-\pi
\tilde{n}_{ij}^2e^{2(y_j-x_i)}\Big\}.\eqa
Computing the integrals  by using
the formula\be \int_{-\infty}^{\infty}e^{\imath\omega
x-px^2}dx=\sqrt{\frac{\pi}{p}}\,\,\,\,e^{-\frac{\omega^2}{4p}}, \qquad
p>0,\ee
we readily obtain
\bqa
\tilde{Q}^{\mathfrak{gl}_{\ell+1}}(\underline{x},\underline{y}|s)=
2^{\ell+1}\,\exp\Big\{\sum_{i=1}^{\ell+1}(\imath s+\rho_i)(y_i-x_i)-\\
\nonumber -\pi
\sum_{i=1}^{\ell}\Big(e^{2(y_i-x_i)}+e^{2(x_{i+1}-y_{i})}\Big)- \pi
e^{2(y_{\ell+1}-x_{\ell+1})}\,\Big\},\eqa where
$\rho_j=\frac{\ell}{2}+1-j,\,\,\,j=1,\ldots,\ell+1$. This completes
the proof of the formula \eqref{twQ} for the integral kernel.

Now we prove \eqref{EigenpropT}.
Taking into account \eqref{Lfactor} we shall check the following
relation:
 \bqa\label{inteqGL}
\prod_{j=1}^{\ell+1}\,
\pi^{-\frac{\i s+\i\lambda_j}{2}}\,
\Gamma\left(\frac{\i s+\imath \lambda_j}{2}\right)
=\int_{G} \,
dg\,\,\tilde{\phi}(g^{-1},s)\,\varphi_{\underline{\lambda}}(g), \eqa
where  ${\rm Im} s <0$ is assumed.
The right  hand side of (\ref{inteqGL})  can be written as follows
 \be
2^{\ell+1}\,\int_{G\times K} \,dk\,dg\,|\det
g|^{\imath s +\frac{\ell}{2}}e^{-\pi {\rm \Tr}
(g^tg)}\,e^{<h(kg),\i \underline{\lambda}-\rho>}=
2^{\ell+1}\,\int_{G} \,dg\,|\det
g|^{\imath s +\frac{\ell}{2}}e^{-\pi {\rm \Tr}
(g^tg)}\,e^{<h(g),\i \underline{\lambda}-\rho>}\\= \nonumber
2^{\ell+1}\,\int_{K\times A\times N_-}\,dn\,d^{\times}a\,dk\,
\delta_{B_-}(a)|\det
a|^{\i s +\frac{\ell}{2}}e^{-\pi\Tr(na^2n^t)}\,
e^{\<\log(a),\i \underline{\lambda}-\rho\>}
\nonumber\ee \be =2^{\ell+1}\,\int_{A\times
N_-}\,dn\,d^{\times}a\,\delta_{B_-}(a)|\det
a|^{\imath s +\frac{\ell}{2}}e^{-\pi\Tr
(na^2n^t)}\,e^{<\log(a),\i \underline{\lambda}-\rho>}\\ =\nonumber
\prod_{j=1}^{\ell+1}\, \pi^{-\frac{\i s -\i\lambda_j}{2}}\,
\Gamma\bigl(\frac{\i s+\i\lambda_j}{2}\bigr),\ee  where the
following formula was used
\be\label{Euler}\int_{-\infty}^{+\infty}dx e^{\nu
x}e^{-ae^{2x}}=\frac{1}{2}\,a^{-\frac{\nu}{2}}\,
\Gamma(\frac{\nu}{2}),\qquad {\rm Re}\, \nu >0,\quad a>0.
\ee
$\Box$

Combining Theorems 2.1, 2.2 one obtains the following result proved in
\cite{PSR}.

\begin{prop}\label{Cor}
Let $\CQ^{(2)}$ be a  $K$-biinvariant function on $GL_{\ell+1}(\IR)$ as
given by 
\be\label{PSRf}
\CQ^{(2)}(g^{-1},s)=\,|\det g|^{\i s+\ell/2}\int_{GL_{\ell+1}(\IR)}
e^{-\pi\Tr Z\,(gg^t\,+1)\,Z^t}\,|\det
Z|^{2\i s+\ell}\,dZ.
\ee
Here $d Z$ is the standard Haar measure on $GL_{\ell+1}(\IR)$.  Then one has the identity
\be\label{Eigenprop2}
\Phi^{\mathfrak{gl}_{\ell+1}}_{\underline{\lambda}}*\CQ^{(2)}(g,s)=
L(s,\underline{\lambda})\,\,
L(s,-\underline{\lambda})\,\,\Phi^{\mathfrak{gl}_{\ell+1}}_{\underline{\lambda}}(g),
\ee
where $L(s,\underline{\lambda})$ is the  local Archimedean $L$-factor
 \eqref{Lfac}.
\end{prop}

\proof Note that the function \eqref{PSRf} can be represented as a
convolution 
$$
\CQ^{(2)}(g,s)=\widetilde{\CQ}^{\mathfrak{gl}_{\ell+1}}*
\CQ^{\mathfrak{gl}_{\ell+1}}(g,s),\qquad
\CQ^{\mathfrak{gl}_{\ell+1}}(g,s)=\,|\det g|^{\i s+\ell/2}\,e^{-\pi\Tr g^tg}, \qquad
\widetilde{\CQ}^{\mathfrak{gl}_{\ell+1}}(g,s)=\CQ^{\mathfrak{gl}_{\ell+1}}(g^{-1},s).
$$ 
The statement of the Proposition directly follows from Theorems \ref{MainTh1},
\ref{MainTh2} (taking into account that 
our normalization of the Haar measure differs from that in \cite{PSR}
by the factor  $2^{\ell+1}$). $\Box$

\section{Universal Baxter operator for classical groups $Sp_{2\ell}$
  and $SO_{2\ell}$ }

In this Section we propose a generalization of Theorem \ref{MainTh1}
to the case of the maximal split forms of the
classical series $SO_{2\ell}$ and $Sp_{2\ell}$. Note that the resulting
expressions for analogs of \eqref{UBO} are not as
simple as in the case of general linear groups  and are given by
non-trivial integrals.

Let us first recall a particular realization of classical Lie groups
$SO_{2\ell}$ and $Sp_{2\ell}$
as subgroups of the general linear group $GL_{2\ell}$.
We are interested in  realizations of  standard representations
$\pi_{2\ell}:\mathfrak{g}\to End(\mathbb{C}^{2\ell})$ of the Lie
algebras $\Fg=\mathfrak{sp}_{2\ell},\,\mathfrak{so}_{2\ell}$  such
that the Weyl generators corresponding to  the Borel (Cartan) subalgebras of
$\Fg$ are realized by upper triangular (diagonal) matrices (see
e.g. \cite{DS}). The corresponding group embedding
$G\to GL_{2\ell}$, $G=SO_{2\ell}\, ,Sp_{2\ell}$ can be defined as follows.
Consider the following
involution on $GL_{2\ell}$ (identified with  a subspace of $2\ell\times
2\ell$-matrices):
\be\label{em1}
g\longmapsto g^*:=S\cdot J\cdot
(g^{-1})^{t}\cdot J^{-1}\cdot S^{-1},
\ee
where $a\to a^{t}$ is induced by the
standard matrix transposition,
\be
S={\rm diag}(1,-1,\ldots,-1,1),
\ee
and  $J=\|J_{i,j}\|=\|\delta_{i+j,2\ell+1}\|$.
The symplectic  group $G=Sp_{2\ell}$ then can be defined as the
following  subgroup of   $GL_{2\ell}$:
\be\label{Spemb}
 Sp_{2\ell}=\{g\in GL_{2\ell}:g^{*}=g\}.
\ee
Similarly in the case of $G=SO_{2\ell}$ consider the
involution on $GL_{2\ell}$
\be\label{em2}
g\longmapsto g^*:=S\cdot J\cdot
(g^{-1})^{t}\cdot J^{-1}\cdot S^{-1},
\ee
where
\be
S={\rm diag}(1,-1,\ldots,(-1)^{\ell-1},(-1)^{\ell-1},(-1)^\ell,\cdots
,1),
\ee
and  $J=\|J_{i,j}\|=\|\delta_{i+j,2\ell+1}\|$.
The  orthogonal  group $G=SO_{2\ell}$  can be defined as the
 following  subgroup of   $GL_{2\ell}$:
\be\label{SOemb}
SO_{2\ell}=\{g\in GL_{2\ell}:g^{*}=g\}.
\ee
The maximal compact subgroup of $G=SO_{2\ell}, Sp_{2\ell}$ embedded
this way is given by an intersection of $G$ with the
maximal compact subgroup of $GL_{2\ell}(\IR)$.

We would like to construct elements of Hecke algebra $\CH(G,K)$,
where $G$ are maximal split forms of $Sp_{2\ell}$, $SO_{2\ell}$ and
$K\subset G$ is a maximal compact subgroup, such that their
actions on spherical vectors in spherical principle series representations are given by
multiplications on the corresponding local Archimedean
$L$-factors associated with the standard representations  of the dual Lie
groups. Recall that $SO_{2\ell}$ is self-dual and
$Sp_{2\ell}$ is dual to $SO_{2\ell+1}$.  The corresponding local $L$-factors are
given by
\be\label{SO}
L^{SO_{2\ell}}(s,\mu_1,r)=\prod_{i=1}^{\ell}
\Gamma_{\IR}(s,\lambda_i)\,
\Gamma_{\IR}(s,-\lambda_i),\qquad G=SO_{2\ell},\quad G^{\vee}=SO_{2\ell},
\ee
\be\label{SP}
L^{SO_{2\ell+1}}(s,\mu_2,r)=
\Gamma_{\IR}(s,0)\,\,\prod_{i=1}^{\ell} \Gamma_{\IR}(s,\lambda_i)
\,\Gamma_{\IR}(s,-\lambda_i),\quad G=Sp_{2\ell},\quad G^{\vee}=SO_{2\ell+1},
\ee
where we parameterize the weights of the corresponding principle series
representations via elements of the dual Lie algebras by taking
$\mu_1={\rm diag}(\lambda_1,\lambda_2,\ldots,
\lambda_{\ell},-\lambda_1,\ldots, -\lambda_{\ell})$ for an element   of
the Lie algebra $\mathfrak{so}_{2\ell}$
and similarly $\mu_2={\rm diag}(0,\lambda_1,\lambda_2,\ldots,
\lambda_{\ell},-\lambda_1,\ldots, -\lambda_{\ell})$
for  an element   of the Lie algebra $\mathfrak{so}_{2\ell+1}$.
Here we  also use  the  modified $\Gamma$-function:
\be\label{modG}
\Gamma_{\IR}(s,\lambda)=\,
\pi^{-\frac{\imath s -\imath\lambda }{2}}
\Gamma\Big(\frac{\imath s-\imath \lambda}{2}\Big).
\ee
The following Theorem is a simple reformulation of the
result obtained in \cite{PSR}.

\begin{te} \label{Main}
Let $G$ be either $SO_{2\ell}$ or $Sp_{2\ell}$ and $G^{\vee}$ be the
corresponding dual group $SO_{2\ell}$ or $SO_{2\ell+1}$.
Let  $\CQ_G(g,s)$ be a one-parameter family of functions on $G$  given by
 \be\label{UBOa}
\CQ_G(g,s)=d_G(s)\,\,\frac{R_G(g,s)}{R_G(0,s)},
\ee
where
\be\label{formula}
R_G(g,s)=\int_{GL_{2\ell}(\IR)}
 e^{-\pi \Tr Z^t( g^tg+1)Z}\,\,|\det Z|^{
   \imath s}\,dZ,
\ee
$$
R_G(0,s)=\int_{GL_{2\ell}(\IR)}
 e^{-\pi\Tr Z^tZ}\,\,|\det Z|^{
  \imath  s}\,dZ,
$$
and
$$
d_{SO_{2\ell}}(s)=\prod^{2\ell-2}_{j\equiv 0({\rm mod} 2),
  j=0}\Gamma_{\IR}(2\imath s-j),
$$
$$
d_{Sp_{2\ell}}(s)=\prod^{2\ell}_{j\equiv 0({\rm mod} 2),
  j=2}\Gamma_{\IR}(2\imath s-j).
$$
We also imply that $G$ is embedded in $GL_{2\ell}$ via
\eqref{Spemb}, \eqref{SOemb}.

Then $\CQ_G(g,s)$ is  a $K$-biinvariant function such that its action
on a spherical
vector $\phi_K$  in a principle series representation
$\CV_{\underline{\mu}}={\rm Ind}_{B_-}^G\,\chi_{\underline{\mu}}$
is via multiplication on the corresponding local
Archimedean $L$-factor (given by \eqref{SO}, \eqref{SP})
\be\label{conS1}
\phi_K*\CQ_G(s)(g)=\int_G\, dg_1\,
\CQ_G(g_1,s)\,\phi_K(g g_1^{-1})=L^{G^{\vee}}(s,\underline{\mu})
\phi_K(g),
\ee
and thus $\CQ_G(s)$ is a universal Baxter operator for classical groups
$SO_{2\ell}$, $Sp_{2\ell}$. In particular it acts on the Whittaker
functions \eqref{equivar}
associated with the classical groups $G$ by a  multiplication on
 the corresponding local Archimedean $L$-factor \eqref{SO},
 \eqref{SP}.
\end{te}

\proof Let $\varphi_{\underline{\mu}}(g)$ be a normalized spherical
function on $G$ corresponding to the principle series representation
$\CV_{\underline{\mu}}$. According to \cite{PSR} (see e.g. p. 37)
the following integral representations
for Archimedean $L$-factors \eqref{SO}, \eqref{SP}
holds
\be\label{PSR}
L^{G^{\vee}}(s,\mu)
=\int_{G}\,\CQ_G(g,s)\, \varphi_{\underline{\mu}}(g^{-1})dg,
\ee
where $\CQ_G(g,s)$ is given by \eqref{UBOa}.
Now using \eqref{conS} and \eqref{Lfactor} we  obtain  \eqref{conS1}. $\Box$

As a simple illustration of the above results let us provide explicit
calculations for $\ell=1$.  The maximal split form of $SO_2$ can be embedded
in $GL_2(\IR)$ via diagonal matrices. Thus the connected component of
the unity can be parametrized  as follows:
$$
g=\begin{pmatrix} e^t & 0 \\ 0 & e^{-t} \end{pmatrix}.
$$
This is in agreement with \eqref{SOemb} with $S={\rm diag}(1,1)$.
The normalized spherical function is given by
$$
\varphi_{\underline{\mu}}(g^{-1}(t))=e^{-\imath \lambda t},
\qquad \underline{\mu}=(\lambda,-\lambda).
$$
Specialization of \eqref{formula} to $G=SO_2$ gives
$$
R_{SO_2}(g,s)=\int_{GL_2(\IR)}\,e^{-\pi\Tr Z^t(g^tg+1)Z}\, |\det
Z|^{\imath s}\,dZ.
$$
To calculate the integral let us denote
$X^2=g^tg+1={\rm diag}(X_1^2,X_2^{2})$ and $X^2_1=e^{2t}+1$,
$X_2^2=e^{-2t}+1$.  Then we have
$$
\Tr Z^tX^2Z=\sum_{i,j=1}^2 Z_{ij}^2X_i^2.
$$
Using the change of variables $Z_{ij}\to Z_{ij}X_i^{-1}$  we
obtain
$$
R_{SO_2)}(g,s)=|X_1X_2|^{-\imath s}\,\int_{GL_2(\IR)}\,e^{-\pi\Tr Z^tZ}\, |\det
Z|^{\imath s}\,dZ.
$$
Thus we have
\be\label{i11}
\CQ_{SO_2}(g,s)=d_{SO_2}(\imath s)\,
\frac{R_{SO_2}(g,s)}{R_{SO_2}(0,s)}=\frac{\Gamma_{\IR}(2\imath s)}
{(e^{2t}+1)(1+e^{-2t})^{\imath s/2}}=\frac{\Gamma_{\IR}(2\imath s)}
{(e^{t}+e^{-t})^{\imath s}}.
\ee
Now we would like to calculate the integral
$$
L^{SO_2}(s,\underline{\mu})=2\int_{\IR}
 \,dt \,e^{- \imath \lambda t}\,\CQ_{SO_2}(g(t),s),
$$
where prefactor $2$ takes into account sum over connected components
of the split form of $SO(2)$.  
Taking into account \eqref{i11} and using the Euler integral representation for
the Gamma-function  we have ( assuming ${\rm Im} s<0$)
$$
L^{SO_2}(s,\underline{\mu})=
2\int_{\IR^2}\,d\tau \,dt \,e^{ \imath s\tau}e^{-\imath \lambda t}
e^{-\pi e^{\tau}(e^t+e^{-t})}\,
$$
$$=\int_{\IR^2} dt_1\,dt_2\,\,
e^{t_1\frac{(\imath s- \imath \lambda)}{2}+t_2\frac{(\imath s+\imath
    \lambda)}{2}} e^{-\pi(e^{t_1}+e^{t_2})}=
\pi^{-\frac{\imath s-\imath \lambda}{2}}
\Gamma\left(\frac{\imath s-\imath \lambda}{2}\right)
\pi^{-\frac{\imath s+\imath \lambda}{2}}\Gamma\left(\frac{\imath s+\imath
  \lambda}{2}\right)
=L_{SO_2}(s,\mu).
$$
One can directly calculate the action on $SO_2$-Whittaker
functions
$$
\Phi_{\underline{\mu}}(x)=e^{\imath \lambda x},\qquad \underline{\mu}=(\lambda,-\lambda).
$$
We have
$$
\Phi_{\underline{\mu}}*\CQ_{SO_2}(s) (x)=\pi^{-\frac{\imath s-\imath \lambda}{2}}
\Gamma\left(\frac{\imath s-\imath \lambda}{2}\right)\pi^{-\frac{\imath s+\imath \lambda}{2}}
\Gamma\left(\frac{\imath s+\imath \lambda}{2}\right)\,\Phi_{\underline{\mu}}(x)=
L^{SO(2)}(s,\underline{\mu})\Phi_{\underline{\mu}}(x).
$$
This completes explicit verification of the statement of Theorem
\ref{Main} for $G=SO_2$.

\vskip 1cm

\noindent {\small {\bf A.G.} {\sl Institute for Theoretical and
Experimental Physics, 117259, Moscow,  Russia; \hspace{8 cm}\,
\hphantom{xxx}  \hspace{2 mm} School of Mathematics, Trinity College
Dublin, Dublin 2, Ireland; \hspace{6 cm}\hspace{5 mm}\,
\hphantom{xxx}   \hspace{2 mm} Hamilton Mathematics Institute,
Trinity College Dublin, Dublin 2, Ireland;}\\

\noindent{\small {\bf D.L.} {\sl
 Institute for Theoretical and Experimental Physics,
117259, Moscow, Russia}.
\end{document}